\documentclass[reqno,11pt]{amsart}
\pdfoutput=1
\usepackage{amsmath,amsthm,amsfonts,amssymb,ifpdf}
\usepackage{amssymb}
\usepackage{amsmath}
\usepackage{amsthm}
\usepackage{graphicx}
\usepackage[all]{xy}
\usepackage{enumerate}
\usepackage{tikz-cd}

\theoremstyle{plain} \newtheorem{introthm}{Theorem}
\newtheorem{theorem}{Theorem}[section]
\newtheorem{proposition}[theorem]{Proposition}
\newtheorem{lemma}[theorem]{Lemma}
\newtheorem{corollary}[theorem]{Corollary}
\newtheorem{question}[theorem]{Question}

\theoremstyle{definition}

\theoremstyle{remark} \newtheorem{remark}[theorem]{Remark}

\usepackage[all]{xy}

\ifpdf
  \usepackage[
    pdftex,
    colorlinks,%
    linkcolor=blue,citecolor=red,urlcolor=blue,
    hyperindex,%
    plainpages=false,%
    bookmarksopen,%
    bookmarksnumbered%
  ]{hyperref} 
 
 \usepackage{thumbpdf}
\else
  \usepackage{hyperref}
\fi

\def\Sym{{\rm Sym}}

\def\SL{{\rm SL}}

\def\dim{{\rm dim}}

\def\id{{\rm id}}

\newcommand{\ncom}{\newcommand}
\ncom{\mylabel}[1]{{\rm (#1)}\label{#1}}
\ncom{\Hom}{{\textit{Hom}}}
\ncom{\eop}{{\hfill $\Box$}}
\begin{document}

\setcounter{tocdepth}{1}

\title[Multiprojective spaces]{A representation theoretic classification of multiprojective spaces}

\author{Arijit Mukherjee}
\address{Department of Mathematics\\ Indian Institute of Technology Madras\\ Tamil Nadu - 600 036, India.}
\email{mukherjee7.arijit@gmail.com}

\begin{abstract}
For a partition $(n_1,\ldots,n_r)$ of a positive integer $n$, consider the associated multiprojective space $\mathbb{P}^{n_1}\times\cdots\times\mathbb{P}^{n_r}$.  That multiprojective spaces attached to distinct partitions of $n$ are pairwise non-isomorphic is known, having been established by algebro-geometric methods.  In this paper we give a new, representation-theoretic proof.  We prove that the complex cohomology ring of a multiprojective space carries a natural $\mathfrak{sl}(2,\mathbb{C})$-module structure compatible with the K\"unneth decomposition. The classification then follows from C. S. Rajan's theorem on the unique decomposition of tensor products of irreducible representations of a simple Lie algebra.  We further show, by a dimension count argument, that this approach is intrinsic to the genus zero curve, distinguishing $\mathbb{P}^1$ from curves of higher genus. 
\end{abstract}
\maketitle


\section{Introduction}
\label{sec:Introduction}
Projective spaces and their products are among the most basic and pivotal objects in algebraic geometry.  They are the most fundamental examples of complete varieties, and they occur throughout the literature.  For instance, in \cite{FM} the authors study the arithmetically Cohen-Macaulay property for sets of points in products of projective spaces, and in \cite{GPP} the authors compute the effective cones of blow ups of certain products of projective spaces in up to six points, showing that the resulting varieties are log Fano.  A product of projective spaces, following the terminology of C. Musili (cf. \cite[Chapter 3, \S 36, p. 150]{Mus}), is called a \emph{multiprojective space}.  A basic question about these objects is the following.
\begin{question}\label{intro question}
Let $r$ and $s$ be positive integers, and let $(n_1,\ldots,n_r)$ and $(m_1,\ldots,m_s)$ be an $r$-tuple and an $s$-tuple of positive integers respectively, distinct from one another.  Are the multiprojective spaces $\mathbb{P}^{n_1}\times \cdots\times\mathbb{P}^{n_r}$ and $\mathbb{P}^{m_1}\times\cdots\times \mathbb{P}^{m_s}$ non-isomorphic?
\end{question}
The answer is affirmative for the trivial reason of dimension if $\sum_{i=1}^r n_i\neq \sum_{j=1}^s m_j$.  We therefore restrict to the equidimensional case, which is the substantive one:
\begin{question}\label{modified question}
Let $n$ be a positive integer.  Given two distinct partitions $(n_1,\ldots,n_r)$ and $(m_1,\ldots,m_s)$ of $n$, are the multiprojective spaces $\mathbb{P}^{n_1}\times \cdots\times\mathbb{P}^{n_r}$ and $\mathbb{P}^{m_1}\times\cdots\times \mathbb{P}^{m_s}$ non-isomorphic?
\end{question}

Here the dimension argument no longer applies.  Question \ref{modified question} has been answered affirmatively, by algebro-geometric means, in \cite{BSV} and \cite{MN}.  In \cite{BSV}, the authors use the Chow rings of the multiprojective spaces, in the course of classifying multiprojective bundles over projective spaces.  In \cite{MN}, the authors use the Picard group and the cohomology ring, interpreting the generators of the cohomology rings as the primitive generators of the nef cones of the corresponding products, to classify multiprojective spaces. Then, using it, they classify products of symmetric products of an arbitrary smooth projective curve over $\mathbb{C}$, while counting the Hilbert schemes of a curve associated to the good partitions of a constant polynomial satisfying the diagonal property.

In this paper we answer Question \ref{modified question} affirmatively by a different, representation-theoretic route.  The strategy rests on three observations.  First, for the projective space the cohomology ring $H^{\ast}(\Sym^n(\mathbb{P}^1),\mathbb{C})=H^{\ast}(\mathbb{P}^n,\mathbb{C})$ carries an $\mathfrak{sl}(2,\mathbb{C})$-module structure, coming from the action of the Lefschetz operators $L$, $\Lambda$ and $h$, and this module is precisely the irreducible module $\Sym^n(\mathbb{C}^2)$ (cf. Lemma \ref{lemma for alternative proof using Rajan}).  Second, the K\"unneth isomorphism is an isomorphism of $\mathfrak{sl}(2,\mathbb{C})$-modules (cf. Lemma \ref{lemma_Kunneth is sl2 equivariant}), so the cohomology of a multiprojective space $\mathbb{P}^{n_1}\times\cdots\times\mathbb{P}^{n_r}$ is the tensor product $\Sym^{n_1}(\mathbb{C}^2)\otimes\cdots\otimes\Sym^{n_r}(\mathbb{C}^2)$ of nontrivial irreducible $\mathfrak{sl}(2,\mathbb{C})$-modules.  Third, an isomorphism of multiprojective spaces induces an isomorphism of these cohomology rings as $\mathfrak{sl}(2,\mathbb{C})$-modules (cf. Corollary \ref{Corollary_isomorphism of cohomology rings as modules}).  By C. S. Rajan's theorem (cf. \cite{Raj}) on the unique decomposition of tensor products of irreducible representations of a simple Lie algebra, the tensor factors, and hence the partition, are determined up to reordering.  This yields our main result (cf. Theorem \ref{Main Theorem_with proof}).
\begin{introthm}\label{intro main theorem}
Let $n$ be a positive integer.  For any two distinct partitions $(n_1,\ldots,n_r)$ and $(m_1,\ldots,m_s)$ of $n$, the multiprojective spaces $\mathbb{P}^{n_1}\times\cdots\times \mathbb{P}^{n_r}$ and $\mathbb{P}^{m_1}\times \cdots\times\mathbb{P}^{m_s}$ are non-isomorphic.
\end{introthm}

This gives a proof of a known result by a method that connects two branches of mathematics, algebraic geometry and representation theory, and that, to the best of our knowledge, has not been used in this context before.

We further delineate the scope of this method.  The $\mathfrak{sl}(2,\mathbb{C})$-module structure above rests on the coincidence $H^{\ast}(\Sym^n(\mathbb{P}^1),\mathbb{C})\cong\Sym^n(H^{\ast}(\mathbb{P}^1,\mathbb{C}))$, which holds for the genus zero curve $\mathbb{P}^1$.  For a smooth projective curve $C$ of genus $g\geq 1$ this coincidence fails.  In fact, we show (cf. Proposition \ref{dimension inequality proposition}) that $\dim_{\mathbb{C}}\;H^{\ast}(\Sym^n(C),\mathbb{C})<\dim_{\mathbb{C}}\;\Sym^n(H^{\ast}(C,\mathbb{C}))$ for all $n\geq 2$.  Therefore, the representation-theoretic argument exploits a feature special to the projective line among all smooth curves.  For higher genus curves, the classical algebro-geometric techniques of \cite{MN} remain the available tool.

\section{Preliminaries}
\subsection{Symmetric product of curve} 
Throughout this article, we take $C$ to be a projective curve of genus $g$ over complex numbers.  For any positive integer $d$, by $C^{\times d}$ we denote the $d$-fold product of the curve $C$.  The notation $\Sym^{d}(C)$ is reserved for the $d$-th symmetric product of $C$.  $\Sym^{d}(C)$ is obtained as the quotient $\frac{C^{\times d}}{\sigma ^d}$ of $C^{\times d}$ under the action of the permutation group $\sigma ^d$ of $d$ symbols.  It is well known that $\Sym^d(C)$ is a smooth projective variety over $\mathbb{C}$ (cf. \cite[Proposition 3.1 and Proposition 3.2, p.~9]{Mi}).  Moreover, it is of dimension $d$ over $\mathbb{C}$ as the quotient map from $C^{\times d}$ to $\Sym^d(C)$ is a finite morphism (cf. \cite[Proposition 3.1, p.~9]{Mi}).

For a genus $0$ curve $C$, i.e., when $C$ is isomorphic to the complex projective line $\mathbb{P}^1$, it is well known that $\Sym^n(\mathbb{P}^1)$ is isomorphic to the $n$-dimensional complex projective space $\mathbb{P}^n$.  Incorporating this fact,  we now recall the cohomology ring structure of $\mathbb{P}^n$ as we use the same several times.
\begin{proposition}\label{cohomology of projective spaces}
The $i$-th graded piece $H^i(\Sym^n(\mathbb{P}^1), \mathbb{C})$ of the cohomology ring $H^{\ast}(\Sym^n\\(\mathbb{P}^1), \mathbb{C})$ of the $n$-th symmetric product of the complex projective line $\mathbb{P}^1$ is given by
\begin{equation*}
H^i(\Sym^n(\mathbb{P}^1), \mathbb{C})=H^i(\mathbb{P}^n, \mathbb{C})= \left\{ \begin{array}{ll}
\mathbb{C},& \mbox{$i$ even, $0\leq i \leq 2n$};\\ 0, &
\mbox{otherwise}.\end{array} \right.
\end{equation*}
\end{proposition}
We now recall the Betti numbers and the generating series of the Poincar\'e polynomial of the symmetric product of a curve of genus greater than zero. 
\begin{proposition}\label{Betti numbers of symmetric product}(\;\cite[Equation 4.2, p.~322]{Mac})
Let $C$ be a smooth projective curve over $\mathbb{C}$ of genus $g\geq 1$.  Then the $r$-th Betti number $B_r$ of $\Sym^n(C)$ is given by 
\begin{equation*}
B_r=B_{2n-r}=\big(\begin{smallmatrix}
  2g\\
  r
\end{smallmatrix}\big)+\big(\begin{smallmatrix}
  2g\\
  r-2
\end{smallmatrix}\big)+\ldots \;\;, 
\end{equation*}
where $0\leq r\leq n$.
\end{proposition}
\begin{proposition}\label{generating series}(\;\cite[Equation 4.3, p.~322]{Mac})
Let $C$ be a smooth projective curve over $\mathbb{C}$ of genus $g\geq 1$.  Then the Poincar\'e polynomial $B_0+B_1x+\cdots +B_{2n}x^{2n}$ of $\Sym^n(C)$ is the coefficient of $t^n$ in the expansion of 
\begin{equation*}
\dfrac{(1+tx)^{2g}}{(1-t)(1-tx^2)}.
\end{equation*}
\end{proposition}
\subsection{Irreducible representations of $\mathfrak{sl}(2,\mathbb{C})$}
Finite dimensional representations of semisimple Lie algebras are very well studied. Here we recall finite dimensional irreducible representations of the Lie algebra $\mathfrak{sl}(2,\mathbb{C})$ of the special linear group $\SL(2,\mathbb{C})$, as these are required in the context of this paper.  The Lie algebra $\mathfrak{sl}(2,\mathbb{C})$ is realised as the $\mathbb{C}$-vector space of $2\times 2$ trace zero matrices with complex entries.  It is equipped with bracket operation given as 
\begin{equation*}
    [A,B]=AB-BA.
\end{equation*}
The standard basis of $\mathfrak{sl}(2,\mathbb{C})$ consists of the following three elements 
\begin{equation}\label{Generators of trace zero matrices}
    X=\begin{pmatrix}
    0 & 1\\
    0 & 0
    \end{pmatrix}, Y=\begin{pmatrix}
    0 & 0\\
    1 & 0
    \end{pmatrix} \text{\;and\;} H=\begin{pmatrix}
    1 & 0\\
    0 & -1
    \end{pmatrix},
\end{equation}
and they satisfy the following bracket relations:
\begin{equation}\label{Bracket relations between matrix generators}
  [X,Y]=H,\;
    [H,X]=2X,\;
    [H,Y]=-2Y.
    \end{equation}
A finite dimensional representation of $\mathfrak{sl}(2,\mathbb{C})$ is nothing but a Lie algebra map $\rho: \mathfrak{sl}(2,\mathbb{C})\rightarrow gl(V) $, where $V$ is a finite dimensional $\mathbb{C}$-vector space and $gl(V)$ is the algebra of endomorphisms of $V$.  Equivalently, $V$ is nothing but a $\mathfrak{sl}(2,\mathbb{C})$-module.  Moreover, if $V$ is irreducible and $\dim_{\mathbb{C}}\;V=n+1$, then $V$ is isomorphic to the direct sum of the one-dimensional eigenspaces corresponding to $n+1$ distinct (integer) eigenvalues of $\rho(H)$ and vice versa.       
\begin{proposition}\label{simplicity of sl(n)}
Let $\mathfrak{sl}(2,\mathbb{C})$ be the Lie algebra of the special linear group $\SL(2,\mathbb{C})$.  Let $\Sym^d(V)$ be the $d$-th symmetric product of a finite dimensional $\mathbb{C}$-vector space $V$.  Then 
\begin{enumerate}
    \item $\mathfrak{sl}(2,\mathbb{C})$ is simple.
    \item There is a bijection between the set of nonnegative integers and the set of irreducible representations of finite dimension, given by 
    \begin{equation*}
        n \mapsto \Sym^n(\mathbb{C}^2).
        \end{equation*}        
\end{enumerate} 
\end{proposition}
\begin{proof}
 See \cite[\S 6.4, p. 125]{Stw} and \cite[Chapter 0, \S 7, p.~118-120]{GH} for $(1)$ and $(2)$ respectively.
\end{proof}

\section{Proof of Main Theorem}

In this section, we prove that any two multi symmetric products of a genus $0$ curve $C$ corresponding to two different partitions of a positive integer are not isomorphic.  We prove this using a classification of tensor products of some modules over some simple Lie algebra, provided by \cite{Raj}.

Before that, we show that this method can not be applicable for curves $C$ of the genus $g>0$, by comparing the symmetric product $\Sym^d(C)$ of a curve $C$ and the symmetric product $\Sym^d(H^{\ast}(C,\mathbb{C}))$ of the vector space $H^{\ast}(C,\mathbb{C})$. 
\begin{proposition}\label{dimension inequality proposition}
Let $C$ be a smooth projective curve over $\mathbb{C}$ of genus $g\geq 1$ and let $n\geq 2$ be an integer.  Then
\begin{equation*}
\dim_{\mathbb{C}} \;H^{\ast}(\Sym^n(C),\mathbb{C})< \dim_{\mathbb{C}}\; \Sym^n(H^{\ast}(C,\mathbb{C})).
\end{equation*}
For $n=1$ the two sides coincide, both being equal to $2g+2$.
\end{proposition}

\begin{proof}
Throughout, we write
\begin{equation*}
L:=\dim_{\mathbb{C}} \;H^{\ast}(\Sym^n(C),\mathbb{C})\qquad\text{and}\qquad R:=\dim_{\mathbb{C}}\; \Sym^n(H^{\ast}(C,\mathbb{C})).
\end{equation*}
We express both $L$ and $R$ as coefficients of $t^n$ in explicit power series and then compare these series coefficientwise.

By Proposition \ref{generating series}, putting $x=1$ in the generating series of the Poincar\'e polynomial of $\Sym^n(C)$ shows that $L$ is the coefficient of $t^n$ in $\tfrac{(1+t)^{2g}}{(1-t)^2}$.  That is,
\begin{equation}\label{L as coefficient}
L=\text{coefficient of\;} t^n \text{\;in\;}\,(1+t)^{2g}(1-t)^{-2}.
\end{equation}
On the other hand, since $\dim_{\mathbb{C}}\;H^{\ast}(C,\mathbb{C})=2g+2$ (cf. \cite[Subsection 2, p.~320]{Mac}), the dimension of the $n$-th symmetric power of a $(2g+2)$-dimensional vector space is
\begin{equation}\label{R as coefficient}
R=\big(\begin{smallmatrix} 2g+n+1\\ n\end{smallmatrix}\big)=\text{coefficient of\;} t^n \text{\;in\;}\,(1-t)^{-(2g+2)},
\end{equation}
where the last equality uses the standard expansion $(1-t)^{-D}=\sum_{k\geq 0}\big(\begin{smallmatrix} D+k-1\\ k\end{smallmatrix}\big)t^k$ with $D=2g+2$.  Factoring the series in \eqref{R as coefficient} as
\begin{equation*}
(1-t)^{-(2g+2)}=(1-t)^{-2}\,(1-t)^{-2g}\;,
\end{equation*}
and comparing this with \eqref{L as coefficient}, both $L$ and $R$ are coefficients of $t^n$ in a product having the common factor $(1-t)^{-2}=\sum_{k\geq 0}(k+1)t^k$, whose coefficients are all positive.  Thus it suffices to compare the coefficients of the remaining factors $(1+t)^{2g}$ and $(1-t)^{-2g}$.  For $j\geq 0$,
\begin{equation*}
\text{coefficient of\;} t^j \text{\;in\;}\,(1+t)^{2g}=\big(\begin{smallmatrix} 2g\\ j\end{smallmatrix}\big)\qquad\text{and}\qquad \text{coefficient of\;} t^j \text{\;in\;}\,(1-t)^{-2g}=\big(\begin{smallmatrix} 2g+j-1\\ j\end{smallmatrix}\big),
\end{equation*}
with the convention $\big(\begin{smallmatrix} 2g\\ j\end{smallmatrix}\big)=0$ for $j>2g$.\\
\textit{Claim}:  We claim that
\begin{equation}\label{binomial comparison}
\big(\begin{smallmatrix} 2g\\ j\end{smallmatrix}\big)\leq \big(\begin{smallmatrix} 2g+j-1\\ j\end{smallmatrix}\big)\qquad\text{for all}\;\; j\geq 0,
\end{equation}
with strict inequality whenever $j\geq 2$.\\ 
\textit{Proof of the claim}: For $j=0$ both sides equal $1$, and for $j=1$ both sides equal $2g$.  For $j\geq 2$, we have:
\begin{equation*}
\big(\begin{smallmatrix} 2g\\ j\end{smallmatrix}\big)=\frac{1}{j!}\prod_{k=0}^{j-1}(2g-k),\qquad \big(\begin{smallmatrix} 2g+j-1\\ j\end{smallmatrix}\big)=\frac{1}{j!}\prod_{k=0}^{j-1}(2g+k).
\end{equation*}
Pairing the factors indexed by the same $k$, we have $2g+k\geq 2g-k$ for every $k$, and $2g+k>2g-k$ for every $k\geq 1$.  Since $j\geq 2$, at least the factor with $k=1$ is strict, so the product on the right strictly exceeds that on the left.  This proves \eqref{binomial comparison}.

Now,  applying the Cauchy product formula with the positive coefficients $(n-j+1)$ of $(1-t)^{-2}$, equations \eqref{L as coefficient} and \eqref{R as coefficient} yield
\begin{equation*}
L=\sum_{j=0}^{n}\big(\begin{smallmatrix} 2g\\ j\end{smallmatrix}\big)(n-j+1)\qquad\text{and}\qquad R=\sum_{j=0}^{n}\big(\begin{smallmatrix} 2g+j-1\\ j\end{smallmatrix}\big)(n-j+1).
\end{equation*}
Subtracting termwise,
\begin{equation}\label{difference formula}
R-L=\sum_{j=0}^{n}\Big(\big(\begin{smallmatrix} 2g+j-1\\ j\end{smallmatrix}\big)-\big(\begin{smallmatrix} 2g\\ j\end{smallmatrix}\big)\Big)(n-j+1).
\end{equation}
By \eqref{binomial comparison},  $\big(\begin{smallmatrix} 2g+j-1\\ j\end{smallmatrix}\big)-\big(\begin{smallmatrix} 2g\\ j\end{smallmatrix}\big)\geq 0$ in \eqref{difference formula} is nonnegative, and $(n-j+1)$ is positive for $0\leq j\leq n$.  Hence $R-L\geq 0$.  For strictness, assume $n\geq 2$ and consider the summand $j=2$, which is present since $2\leq n$.  Now, as
\begin{equation*}
\big(\begin{smallmatrix} 2g+1\\ 2\end{smallmatrix}\big)-\big(\begin{smallmatrix} 2g\\ 2\end{smallmatrix}\big)=\big(\begin{smallmatrix} 2g\\ 1\end{smallmatrix}\big)=2g\geq 2>0,
\end{equation*}
using Pascal's rule $\big(\begin{smallmatrix} 2g+1\\ 2\end{smallmatrix}\big)=\big(\begin{smallmatrix} 2g\\ 2\end{smallmatrix}\big)+\big(\begin{smallmatrix} 2g\\ 1\end{smallmatrix}\big)$, and $(n-1)\geq 1$, the $j=2$ summand of \eqref{difference formula} is positive.  Therefore $R-L>0$, which is precisely the desired strict inequality.

For the case $n=1$,  as $\Sym^1(C)=C$, so $L=\dim_{\mathbb{C}}\;H^{\ast}(C,\mathbb{C})=2g+2$.  As $\Sym^1(H^{\ast}(C,\mathbb{C}))$ is nothing but $H^{\ast}(C,\mathbb{C})$ as vector spaces, so $R=2g+2$.  Thus, the two sides coincide for $n=1$, and \eqref{difference formula} indeed gives $R-L=0$ in this case.
\end{proof}

In contrast to the higher genus case just discussed, something stronger holds for genus zero curves.  We encode these properties exclusively enjoyed by genus zero curves through a lemma.  Before that, we introduce some necessary definitions and notations.

Let $A^{\ast}(\mathbb{P}^n)$ be the space of differential forms on $\mathbb{P}^n$ and $A^r(\mathbb{P}^n)$ be its $r$-th graded piece.  Then we have
$$A^r(\mathbb{P}^n)=\oplus_{p+q=r}A^{p,q}(\mathbb{P}^n),$$
where $A^{p,q}(\mathbb{P}^n)$ is the space of differential forms on $\mathbb{P}^n$ of type $(p,q)$.  Let $\Pi^{p,q}: A^{\ast}(\mathbb{P}^n)\rightarrow A^{p,q}(\mathbb{P}^n)$ be the usual projection. Let $\Pi^r$ be the following map :
\begin{equation}\label{Definition of direct sum of projection maps}
    \Pi^r:=\oplus_{p+q=r} \Pi^{p,q}: A^{\ast}(\mathbb{P}^n)\rightarrow A^r(\mathbb{P}^n). 
\end{equation}

We now recall the definition of two operators on $A^{p,q}(\mathbb{P}^n)$, denoted by $L$ and $\Lambda$, as these two operators are useful to make $H^{\ast}(\mathbb{P}^n, \mathbb{C})$ into a $\mathfrak{sl}(2, \mathbb{C})$-module.
\begin{equation}\label{Definition of L}
\begin{split}
L : A^{p,q}(\mathbb{P}^n) & \rightarrow A^{p+1,q+1}(\mathbb{P}^n)\\
\eta &\mapsto \eta \wedge \omega,  
\end{split}
\end{equation}
where $\omega$ is the $(1,1)$-form associated with the Fubini-Study metric on $\mathbb{P}^n$, and 
\begin{equation}\label{Definition of Lambda}
\Lambda (:=L^{\ast}) : A^{p,q}(\mathbb{P}^n)  \rightarrow A^{p-1,q-1}(\mathbb{P}^n)
\end{equation}
be the adjoint operator of $L$.  Let $h$ be the operator defined as follows:
\begin{equation}\label{Definition of h}
    h=\sum_{p=0}^{2n} (n-p)\Pi^p,
\end{equation}
where $\Pi^p$ is as in \eqref{Definition of direct sum of projection maps}.  We have the following relations between the operators $L$, $\Lambda$  and $h$:
\begin{equation}\label{Bracket relations between the operators}
    [\Lambda,L]=h,\;
    [h,L]=-2L,\;
    [h, \Lambda]=2\Lambda.
\end{equation}
Let $\Delta_d(=dd^{\ast}+d^{\ast}d)$ be the $d$-Laplacian and $\mathcal{H}^{\ast}_d(\mathbb{P}^n)$ be the following harmonic space:
\begin{equation*}
\mathcal{H}^{\ast}_d(\mathbb{P}^n)=\{\eta \in A^{\ast} (\mathbb{P}^n)\mid \Delta_d(\eta)=0\}.    
\end{equation*}
As the Lefschetz operators $L$, $\Lambda$ and $h$ commute with $\Delta_d$, they act on $\mathcal{H}^{\ast}_d(\mathbb{P}^n)$ with relations as in \eqref{Bracket relations between the operators}.  As a result of this action, $H^{\ast}(\Sym^n(\mathbb{P}^1),\mathbb{C})$ becomes a $\mathfrak{sl}(2,\mathbb{C})$-module.   
\begin{lemma}\label{lemma for alternative proof using Rajan}
Let $n$ be any positive integer.  Then
\begin{enumerate}
    \item $H^{\ast}(\Sym^n(\mathbb{P}^1),\mathbb{C})\cong\Sym^n(H^{\ast}(\mathbb{P}^1,\mathbb{C}))$, as $\mathbb{C}$-vector spaces.
    \item $H^{\ast}(\Sym^n(\mathbb{P}^1),\mathbb{C})\cong\Sym^n(H^{\ast}(\mathbb{P}^1,\mathbb{C}))$, as $\mathfrak{sl}(2,\mathbb{C})$-modules.  Hence, $H^{\ast}(\Sym^n(\mathbb{P}^1),\mathbb{C})$ is an irreducible $\mathfrak{sl}(2,\mathbb{C})$-module. 
\end{enumerate}
\end{lemma}
\begin{proof}
(1) From Proposition \ref{cohomology of projective spaces}, we have 
\begin{equation*}
H^{\ast}(\mathbb{P}^1, \mathbb{C})=\bigoplus_{i=0}^1H^{2i}(\mathbb{P}^1, \mathbb{C})=\mathbb{C}\oplus \mathbb{C}.
\end{equation*}
Therefore,
\begin{equation}\label{dim_match_genus zero_1}
\dim_{\mathbb{C}}\;\Sym^n(H^{\ast}(\mathbb{P}^1, \mathbb{C}))=\big(\begin{smallmatrix}
  2+n-1\\
  n
\end{smallmatrix}\big)=n+1
\end{equation}
Again from Proposition \ref{cohomology of projective spaces} we have,
\begin{equation*}
\begin{split}
H^{\ast}(\Sym^n(\mathbb{P}^1), \mathbb{C})=H^{\ast}(\mathbb{P}^n, \mathbb{C})&=\bigoplus_{i=0}^nH^{2i}(\mathbb{P}^n, \mathbb{C})\\
&=\underbrace{\mathbb{C}\oplus \mathbb{C}\oplus \cdots \oplus \mathbb{C}}_{n+1\; \text{times}}.
\end{split}
\end{equation*}
Therefore,
\begin{equation}\label{dim_match_genus zero_2}
\dim_{\mathbb{C}}\;H^{\ast}(\Sym^n(\mathbb{P}^1), \mathbb{C})=n+1
\end{equation}
Therefore, by \eqref{dim_match_genus zero_1} and \eqref{dim_match_genus zero_2}, we have the assertion.\\
(2) We have $\mathcal{H}^{\ast}_d(\mathbb{P}^n)\cong H^{\ast}(\mathbb{P}^n, \mathbb{C})$ (cf. \cite[p.~121]{GH}).  Therefore, as \eqref{Bracket relations between matrix generators} and \eqref{Bracket relations between the operators} hold, $H^{\ast}(\mathbb{P}^n, \mathbb{C})(=H^{\ast}(\Sym^n(\mathbb{P}^1), \mathbb{C}))$ becomes a $\mathfrak{sl}(2,\mathbb{C})$-module by sending
\begin{equation*}
    X \rightarrow \Lambda, Y \rightarrow L \text{\;and\;} H \rightarrow h,
\end{equation*}
where $X$, $Y$ and $H$ are as in \eqref{Generators of trace zero matrices} and $L$, $\Lambda$, $h$ are as defined in \eqref{Definition of L}, \eqref{Definition of Lambda}, \eqref{Definition of h} respectively.  Moreover, the eigenspace for $h$ with respect to the eigenvalue $(n-p)$ is $H^{p}(\Sym^n(\mathbb{P}^1), \mathbb{C})$ (cf. \cite[p.~122]{GH}).  Now, by Proposition \ref{cohomology of projective spaces}, the space $H^{p}(\Sym^n(\mathbb{P}^1),\mathbb{C})$ is one-dimensional for $p=0,2,4,\ldots,2n$ and is zero otherwise.  Hence, as an $h$-module, $H^{\ast}(\Sym^n(\mathbb{P}^1),\mathbb{C})$ has precisely the weights $n,n-2,\ldots,-n$, each occurring with multiplicity one.  By the classification of finite dimensional $\mathfrak{sl}(2,\mathbb{C})$-modules (cf. \cite[Chapter 0, \S 7, p.~119-120]{GH}), the irreducible module of dimension $n+1$ is the unique finite dimensional $\mathfrak{sl}(2,\mathbb{C})$-module whose $h$-weights are exactly $n,n-2,\ldots,-n$, each with multiplicity one. Therefore, $H^{\ast}(\Sym^n(\mathbb{P}^1),\mathbb{C})$ is an irreducible $\mathfrak{sl}(2,\mathbb{C})$-module.  Therefore, as $H^{\ast}(\Sym^n(\mathbb{P}^1),\mathbb{C})\cong\Sym^n(H^{\ast}(\mathbb{P}^1,\mathbb{C}))$ as $\mathbb{C}$-vector spaces, the assertion then follows from Proposition \ref{simplicity of sl(n)}.  
\end{proof}
We set up some notations.  Let $N$ be a compact K\"ahler manifold and $\omega_N$ be a K\"ahler form on $N$.  Then by $L_N$, we denote the following operator 
\begin{equation*}\label{definition of L_N}
\begin{split}
L_N: A^{p,q}(N)&\rightarrow A^{p+1,q+1}(N)\\
    \eta &\mapsto \eta \wedge \omega_N.
\end{split}
\end{equation*}
By $\Lambda_N$ and $h_N$, we denote the operators $L^{\ast}_N$ and $[\Lambda_N,L_N]$ respectively.  Now, if $M$ is another compact K\"ahler manifold and  $\phi:M\rightarrow N$ is an isomorphism, then $\omega_M:=\phi^{\ast}\omega_N$ is again a K\"ahler form on $M$ as $\phi$ is a biholomorphism.  In that case, by $L_M$, $\Lambda_M$ and $h_M$ we denote corresponding operators on $M$ with respect to the K\"ahler form $\omega_M$. 
\begin{lemma}\label{Lemma_functoriality of Lefschetz operators}
Let $M$ and $N$ be two compact K\"ahler manifolds and $\phi:M\rightarrow N$ an isomorphism, and equip $M$ with the K\"ahler form $\omega_M:=\phi^{\ast}\omega_N$. Then the following holds :
 \begin{enumerate}
     \item $\phi^{*}L_N=L_M\phi^{*}$.
     \item $\phi^{*}\Lambda_N=\Lambda_M\phi^{*}$.
     \item $\phi^{*}h_N=h_M\phi^{*}$.
 \end{enumerate}
\end{lemma}
\begin{proof}
(1)   For any $\alpha\in H^{*}(N,\mathbb{C})$,
\begin{align*}
\phi^{*}(L_N\alpha) &= \phi^{*}( \alpha \wedge \omega_N)= \phi^{*}\alpha \wedge \phi^{*}\omega_N = \phi^{*}\alpha \wedge \omega_M = L_M(\phi^{*}\alpha).
\end{align*}
Therefore, $\phi^{*}L_N=L_M\phi^{*}$.\\
(2) Since $\phi:M\rightarrow N$ is a biholomorphism and the K\"ahler form on $M$ is taken to be $\omega_M=\phi^{\ast}\omega_N$, the map $\phi$ is an isometry of the underlying K\"ahler metrics.  Consequently $\phi^{\ast}$ commutes with the Hodge star operators (cf. \cite[Chapter 0, \S 6, p.~82]{GH}), $\phi^{\ast}\circ *_N=*_M\circ \phi^{\ast}$, and hence preserves the Hodge inner products (cf. \cite[Chapter 0, \S 6, p.~80]{GH}). That is, denoting by $\langle\cdot,\cdot\rangle_M$ and $\langle\cdot,\cdot\rangle_N$ the Hodge inner products on $M$ and $N$ respectively, we have
\begin{equation*}
\langle\phi^{\ast}\beta,\phi^{\ast}\gamma\rangle_M=\langle\beta,\gamma\rangle_N.
\end{equation*}
The operator $\Lambda_N$ is, by definition, the adjoint $L_N^{\ast}$ of $L_N$ with respect to $\langle\cdot,\cdot\rangle_N$, and likewise $\Lambda_M=L_M^{\ast}$ with respect to $\langle\cdot,\cdot\rangle_M$.  Therefore, for all $\beta,\gamma\in H^{\ast}(N,\mathbb{C})$,
\begin{align*}
\langle\phi^{*}(\Lambda_N\beta),\phi^{*}\gamma\rangle_M &= \langle\Lambda_N\beta,\gamma\rangle_N = \langle\beta,L_N\gamma\rangle_N\\
&= \langle\phi^{*}\beta,\phi^{*}(L_N\gamma)\rangle_M\\
&= \langle\phi^{*}\beta,L_M\phi^{*}\gamma\rangle_M \qquad (\text{by\; part\;} (1))\\
&= \langle\Lambda_M\phi^{*}\beta,\phi^{*}\gamma\rangle_M.
\end{align*}
As $\phi^{\ast}$ is an isomorphism, $\phi^{\ast}\gamma$ ranges over all of $H^{\ast}(M,\mathbb{C})$, so nondegeneracy of $\langle\cdot,\cdot\rangle_M$ implies $\phi^{*}(\Lambda_N\beta)=\Lambda_M\phi^{*}\beta$ for every $\beta$.  Hence $\phi^{*}\Lambda_N=\Lambda_M\phi^{*}$.\\
(3) Finally, part (1) and part (2), together with the relations $h_M=[\Lambda_M, L_M]$ and $h_N=[\Lambda_N, L_N]$ give $\phi^{*}h_N=h_M\phi^{*}$.
\end{proof}
The following result says that if two multiprojective spaces are isomorphic, then their complex cohomology rings are not only isomorphic as $\mathbb{C}$-vector spaces but also isomorphic as $\mathfrak{sl}(2,\mathbb{C})$-modules.
\begin{corollary}\label{Corollary_isomorphism of cohomology rings as modules}
Let the multiprojective spaces $\mathbb{P}^{n_1}\times \cdots\times\mathbb{P}^{n_r}$ and $\mathbb{P}^{m_1}\times\cdots\times \mathbb{P}^{m_s}$ be isomorphic.  Then, the corresponding cohomology rings $H^{\ast}(\mathbb{P}^{n_1}\times\cdots\times \mathbb{P}^{n_r}, \mathbb{C})$ and $H^{\ast}(\mathbb{P}^{m_1}\times\cdots\times \mathbb{P}^{m_s}, \mathbb{C})$ are isomorphic as $\mathfrak{sl}(2,\mathbb{C})$-modules.
\end{corollary}
\begin{proof}
    Let $\phi$ be an isomorphism between the multiprojective spaces $M:=\mathbb{P}^{n_1}\times \cdots\times\mathbb{P}^{n_r}$ and $N:=\mathbb{P}^{m_1}\times\cdots\times \mathbb{P}^{m_s}$.  Then, clearly $H^{\ast}(M,\mathbb{C})$ and $H^{\ast}(N,\mathbb{C})$ are isomorphic as $\mathbb{C}$-vector spaces for dimension reasons.  It can be noted further that $\mathfrak{sl}(2,\mathbb{C})$ acts on $H^{\ast}(M,\mathbb{C})$ and $H^{\ast}(N,\mathbb{C})$ exactly the same way as mentioned in the proof of the second part of Lemma \ref{lemma for alternative proof using Rajan}.  That is to say, $H^{\ast}(M,\mathbb{C})$ and $H^{\ast}(N,\mathbb{C})$ become $\mathfrak{sl}(2,\mathbb{C})$-modules by sending
    \begin{equation*}
    \begin{split}
    X &\rightarrow \Lambda_M, Y \rightarrow L_M, H \rightarrow h_M,\\
    X &\rightarrow \Lambda_N, Y \rightarrow L_N, H \rightarrow h_N,
    \end{split}
\end{equation*}
 where $X$, $Y$ and $H$ are as in \eqref{Generators of trace zero matrices}, and $L_N$, $\Lambda_N$, $h_N$, $L_M$, $\Lambda_M$, $h_M$ are as defined just before Lemma \ref{Lemma_functoriality of Lefschetz operators}.   Therefore, as $\{X,Y,H\}$ is a basis of $\mathfrak{sl}(2,\mathbb{C})$, by Lemma \ref{Lemma_functoriality of Lefschetz operators}, we have :
\begin{equation*}
\phi^{*}(A\cdot\alpha)=A\cdot(\phi^{*}\alpha),
\end{equation*}
for all $A\in\mathfrak{sl}(2,\mathbb{C})$ and $\alpha\in H^{*}(N,\mathbb{C})$.  That is, the isomorphism $\phi$ preserves the $\mathfrak{sl}(2,\mathbb{C})$ action as well.  Hence, the assertion follows.
\end{proof}
The following result is quite interesting as it completely classifies the tensor products of nontrivial, irreducible, finite dimensional modules over any simple Lie algebra over $\mathbb{C}$. 
\begin{proposition}\label{Rajan result on unique decomposition}
Let $\mathfrak{g}$ be a simple Lie algebra over complex numbers.  Let $V_i$'s and $W_j$'s be nontrivial, irreducible, finite dimensional $\mathfrak{g}$-modules, where $1\leq i \leq m$, $1\leq j\leq n$.  Then
\begin{equation*}
\begin{split}
V_1&\otimes \cdots \otimes V_m \simeq W_1\otimes \cdots \otimes W_n \text{\;as\;} \mathfrak{g}\text{-modules\;}\\
&\Rightarrow m=n \;\& \text{\;there\;exists\;some\;} \sigma \in S_n \text{\;such\;that\;} V_i\simeq W_{\sigma(i)} \text{\;as\;} \mathfrak{g}\text{-modules\;}.
\end{split}
\end{equation*}
\end{proposition}
\begin{proof}
See \cite{Raj}.
\end{proof}

Let $M_1,\ldots,M_r$ be compact K\"ahler manifolds, $M:=M_1\times\cdots\times M_r$, and $\pi_i:M\rightarrow M_i$ the projections.  We now check that the K\"unneth isomorphism $H^{\ast}(M,\mathbb{C})\cong H^{\ast}(M_1,\mathbb{C})\otimes\cdots\otimes H^{\ast}(M_r,\mathbb{C})$
of cohomology rings of compact K\"ahler manifolds is in fact an $\mathfrak{sl}(2,\mathbb{C})$-module isomorphism.  To do so, we first note an important observation.  Clearly, each summand $\pi_i^{\ast}\omega_{M_i}$ is a closed form of type $(1,1)$, and is positive along the $i$-th factor and vanishes along the others.  However, the sum $\omega_M=\sum_{i=1}^{r}\pi_i^{\ast}\omega_{M_i}$ is positive.  Indeed, at any point $x=(x_1,\ldots,x_r)\in M$ the holomorphic tangent space splits as $T_xM=\bigoplus_{i=1}^{r}T_{x_i}M_i$, so a nonzero tangent vector $v=(v_1,\ldots,v_r)$ has $v_k\neq 0$ for some $k$, and
\begin{equation*}
\omega_M(v,\bar v)=\sum_{i=1}^{r}\omega_{M_i}(v_i,\bar v_i)\geq \omega_{M_k}(v_k,\bar v_k)>0\;.
\end{equation*}
Therefore, $\omega_M$ is  a K\"ahler form on $M$.
\begin{lemma}\label{lemma_Kunneth is sl2 equivariant}
Let $M_1,\ldots,M_r$ be compact K\"ahler manifolds and let $M:=M_1\times\cdots\times M_r$ be equipped with the product K\"ahler form $\omega_M:=\sum_{i=1}^{r}\pi_i^{\ast}\omega_{M_i}$, where $\pi_i:M\rightarrow M_i$ denotes the $i$-th projection.  Then the K\"unneth isomorphism
\begin{equation*}
H^{\ast}(M,\mathbb{C})\cong H^{\ast}(M_1,\mathbb{C})\otimes\cdots\otimes H^{\ast}(M_r,\mathbb{C})
\end{equation*}
is an isomorphism of $\mathfrak{sl}(2,\mathbb{C})$-modules, where the left hand side carries the Lefschetz $\mathfrak{sl}(2,\mathbb{C})$-action associated with $\omega_M$ and the right hand side carries the tensor product of the Lefschetz $\mathfrak{sl}(2,\mathbb{C})$-actions of the factors.
\end{lemma}
\begin{proof}
We write $L_i$, $\Lambda_i$, $h_i$ for the Lefschetz operators on $H^{\ast}(M_i,\mathbb{C})$ and $L$, $\Lambda$, $h$ for those on $H^{\ast}(M,\mathbb{C})$.  For $1\leq i\leq r$, set
\begin{equation*}
L^{(i)}:=\id\otimes\cdots\otimes L_i\otimes\cdots\otimes\id,\qquad \Lambda^{(i)}:=\id\otimes\cdots\otimes \Lambda_i\otimes\cdots\otimes\id,\qquad h^{(i)}:=\id\otimes\cdots\otimes h_i\otimes\cdots\otimes\id,
\end{equation*}
where the nontrivial factor sits in the $i$-th slot.  Recall that $L_i$ raises the total degree by $2$, $\Lambda_i$ lowers it by $2$, and $h_i$ preserves it.  In particular, all three are operators of even degree.  By the K\"unneth isomorphism, the classes of the form $\alpha_1\otimes\cdots\otimes\alpha_r$ , with $\alpha_i\in H^{\ast}(M_i,\mathbb{C})$, span $H^{\ast}(M,\mathbb{C})$.  Since all the operators appearing below are $\mathbb{C}$-linear, it suffices to verify the operator identities $L=\sum_i L^{(i)}$, $\Lambda=\sum_i\Lambda^{(i)}$ and $h=\sum_i h^{(i)}$ on such classes.\\
\textit{Operator identity for $L$}:  Under the K\"unneth identification, a class $\alpha_1\otimes\cdots\otimes\alpha_r$ corresponds to $\pi_1^{\ast}\alpha_1\wedge\cdots\wedge\pi_r^{\ast}\alpha_r$.  Since $\omega_{M_i}$ has degree $2$, the class $\pi_i^{\ast}\omega_{M_i}$ is of even degree and hence commutes with each $\pi_j^{\ast}\alpha_j$ in $H^{\ast}(M,\mathbb{C})$.  Therefore,
\begin{align*}
L(\alpha_1\otimes\cdots\otimes\alpha_r)&=\Big(\sum_{i=1}^{r}\pi_i^{\ast}\omega_{M_i}\Big)\wedge(\pi_1^{\ast}\alpha_1\wedge\cdots\wedge\pi_r^{\ast}\alpha_r)\\
&=\sum_{i=1}^{r}\pi_1^{\ast}\alpha_1\wedge\cdots\wedge\pi_i^{\ast}(\omega_{M_i}\wedge\alpha_i)\wedge\cdots\wedge\pi_r^{\ast}\alpha_r\\
&=\sum_{i=1}^{r}\alpha_1\otimes\cdots\otimes L_i\alpha_i\otimes\cdots\otimes\alpha_r=\sum_{i=1}^{r}L^{(i)}(\alpha_1\otimes\cdots\otimes\alpha_r).
\end{align*}
\textit{Operator identity for $\Lambda$}:  The space of harmonic forms on $M$ is the tensor product of the spaces of harmonic forms on the factors, that is, $\mathcal{H}_d^{\ast}(M)\cong\mathcal{H}_d^{\ast}(M_1)\otimes\cdots\otimes\mathcal{H}_d^{\ast}(M_r)$.  Under this identification, the Hodge inner product on $M$ is the tensor product of the Hodge inner products on the factors.  That is,
\begin{equation*}
\langle\alpha_1\otimes\cdots\otimes\alpha_r,\;\beta_1\otimes\cdots\otimes\beta_r\rangle_M=\prod_{i=1}^{r}\langle\alpha_i,\beta_i\rangle_{M_i}.
\end{equation*}
 Taking adjoints with respect to this inner product, the adjoint of $L^{(i)}=\id\otimes\cdots\otimes L_i\otimes\cdots\otimes\id$ is $\id\otimes\cdots\otimes L_i^{\ast}\otimes\cdots\otimes\id=\Lambda^{(i)}$, since the identity operator on each other slot is self-adjoint.  Hence,
\begin{equation*}
\Lambda=L^{\ast}=\Big(\sum_{i=1}^{r}L^{(i)}\Big)^{\ast}=\sum_{i=1}^{r}\Lambda^{(i)}.
\end{equation*}
\textit{Operator identity for $h$}:  We have the following equality:
\begin{equation*}
h=[\Lambda,L]=\Big[\sum_{j=1}^{r}\Lambda^{(j)},\;\sum_{i=1}^{r}L^{(i)}\Big]=\sum_{i,j=1}^{r}\big[\Lambda^{(j)},L^{(i)}\big].
\end{equation*}
For $i\neq j$, the operators $L^{(i)}$ and $\Lambda^{(j)}$ act nontrivially on different tensor slots and are both of even degree.  Therefore, they commute with each other and thus $[\Lambda^{(j)},L^{(i)}]=0$ for $i\neq j$.  For $i=j$,
\begin{equation*}
[\Lambda^{(i)},L^{(i)}]=\id\otimes\cdots\otimes[\Lambda_i,L_i]\otimes\cdots\otimes\id=\id\otimes\cdots\otimes h_i\otimes\cdots\otimes\id=h^{(i)}.
\end{equation*}
Therefore only the diagonal terms survive and
\begin{equation*}
h=\sum_{i=1}^{r}h^{(i)}.
\end{equation*}
Under the correspondence $X\rightarrow\Lambda$, $Y\rightarrow L$, $H\rightarrow h$, the three operator identities
\begin{equation*}
L=\sum_{i=1}^{r}L^{(i)},\qquad \Lambda=\sum_{i=1}^{r}\Lambda^{(i)},\qquad h=\sum_{i=1}^{r}h^{(i)}
\end{equation*}
say exactly that $X$, $Y$ and $H$ act on $H^{\ast}(M,\mathbb{C})$ as they do in the tensor product $\mathfrak{sl}(2,\mathbb{C})$-module structure on $H^{\ast}(M_1,\mathbb{C})\otimes\cdots\otimes H^{\ast}(M_r,\mathbb{C})$.  As $\{X,Y,H\}$ is a basis of $\mathfrak{sl}(2,\mathbb{C})$, the K\"unneth isomorphism intertwines the two actions, and is therefore an isomorphism of $\mathfrak{sl}(2,\mathbb{C})$-modules.
\end{proof}

Now we are in a stage to provide a classification of the multi symmetric products corresponding to distinct partitions for a genus zero curve through the lens of representation theory.
\begin{theorem}\label{Main Theorem_with proof}
Let $n$ be any positive integer.  Given any two distinct partitions $(n_1,\ldots,n_r)$ and $(m_1,\ldots,m_s)$ of $n$, corresponding multiprojective spaces $\mathbb{P}^{n_1}\times\cdots\times \mathbb{P}^{n_r}$ and $\mathbb{P}^{m_1}\times \cdots\times\mathbb{P}^{m_s}$ are non-isomorphic.
\end{theorem}
\begin{proof}
We denote $H^{\ast}(\mathbb{P}^1,\mathbb{C})$ by $V$ and by $\mathfrak{sl}(V)$ we simply mean $\mathfrak{sl}(2,\mathbb{C})$ as $\dim_{\mathbb{C}}\;V=2$.  Then by Lemma \ref{lemma for alternative proof using Rajan} and Lemma \ref{lemma_Kunneth is sl2 equivariant}, we have the following isomorphisms of $\mathfrak{sl}(2,\mathbb{C})$-modules:
\begin{equation}\label{equation 1_alternative proof using Rajan}
\begin{split}
H^{\ast}(\Sym^{n_1}(\mathbb{P}^1)\times \cdots \times \Sym^{n_r}(\mathbb{P}^1),\mathbb{C})&=H^{\ast}(\Sym^{n_1}(\mathbb{P}^1),\mathbb{C})\otimes \cdots \otimes H^{\ast}(\Sym^{n_r}(\mathbb{P}^1,\mathbb{C}))\\
&=\Sym^{n_1}(V)\otimes \cdots \otimes \Sym^{n_r}(V),
\end{split}
\end{equation}
where the product $\prod_{i=1}^{r}\Sym^{n_i}(\mathbb{P}^1)=\prod_{i=1}^{r}\mathbb{P}^{n_i}$ is equipped with the product K\"ahler form $\sum_{i=1}^{r}\pi_i^{\ast}\omega_i$, $\omega_i$ being the Fubini-Study form on the $i$-th factor.  Also, each $\Sym^{n_i}(V)$ is an irreducible  $\mathfrak{sl}(V)$-module by Proposition \ref{simplicity of sl(n)}.  Moreover, these are nontrivial as $\dim_{\mathbb{C}}\;\Sym^{n_i}(V)=n_i+1\geq 2$. 

Similarly, we have the following isomorphisms of $\mathfrak{sl}(2,\mathbb{C})$-modules:
\begin{equation}\label{equation 2_alternative proof using Rajan}
\begin{split}
 H^{\ast}(\Sym^{m_1}(\mathbb{P}^1)\times \cdots \times \Sym^{m_s}(\mathbb{P}^1),\mathbb{C})&=H^{\ast}(\Sym^{m_1}(\mathbb{P}^1),\mathbb{C})\otimes \cdots \otimes H^{\ast}(\Sym^{m_s}(\mathbb{P}^1,\mathbb{C}))\\
 &=\Sym^{m_1}(V)\otimes \cdots \otimes \Sym^{m_s}(V),
 \end{split}
\end{equation}
where each $\Sym^{m_j}(V)$ is nontrivial irreducible $\mathfrak{sl}(V)$-module. 

Proceeding contrapositively, the proof now follows from Corollary \ref{Corollary_isomorphism of cohomology rings as modules}, Proposition \ref{Rajan result on unique decomposition}, Proposition \ref{simplicity of sl(n)}, \eqref{equation 1_alternative proof using Rajan} \& \eqref{equation 2_alternative proof using Rajan}. 
\end{proof}
\begin{remark}\label{concluding remark}
    Proposition \ref{dimension inequality proposition} says that Lemma \ref{lemma for alternative proof using Rajan} does not hold for higher genus curves.  As a result, the technique used in Theorem \ref{Main Theorem_with proof}  is exclusive to the genus zero curve.  In fact, standard classical techniques come handy to take care of curves of arbitrary genus (cf. \cite[Proposition 4.20, p.~20]{MN}).
\end{remark}
\section*{Acknowledgements}
The author would like to express his gratitude to Prof. D. S. Nagaraj for many useful suggestions throughout the preparation of the manuscript.  The author also thanks Prof. Amritanshu Prasad for clearing a couple of doubts.  The author would like to acknowledge Indian Institute of Science Education and Research Tirupati (Award No. - IISER-T/Offer/PDRF/A.M./M/01/2021) and Indian Institute of Technology Madras (Office order No.F.ARU/R10/IPDF/2024) for financial support.


\begin{thebibliography}{01234567}
\bibitem[1]{BSV} Bansal, A., Sarkar, S. and Vats, S. : {\it Isomorphism of Multiprojective Bundles and Projective Towers}, (2023). (\url{https://arxiv.org/abs/2311.00999}) 

\bibitem[2]{FM} Favacchio, G. and Migliore, J. :  {\it Multiprojective spaces and the arithmetically Cohen-Macaulay property}, Math. Proc. Cambridge Philos. Soc. \textbf{166} (2019), no. 3, 583--597. (\url{https://doi.org/10.1017/S0305004118000142})

\bibitem[3]{GPP} Grange, Tim. Postinghel, E. and Prendergast-Smith, A. : {\it Log Fano blowups of mixed products of projective spaces and their effective cones}, Rev. Mat. Complut. \textbf{36} (2023), no. 2, 393--442.
(\url{https://link.springer.com/article/10.1007/s13163-022-00425-2})

\bibitem[4]{GH}  Griffiths, P. and Harris, J. :  {\it Principles of Algebraic Geometry}, Wiley Classics Library, John Wiley \& Sons. Inc., New York, (1994). (\url{https://onlinelibrary.wiley.com/doi/book/10.1002/9781118032527})

\bibitem[5]{Mac} Macdonald, I. G. : {\it Symmetric products of an algebraic curve}, Topology 1, (1962)  319--343. (\url{https://doi.org/10.1016/0040-9383(62)90019-8})

\bibitem[6]{Mi} Milne, J. :  {\it Jacobian varieties}, (2018). (\url{https://www.jmilne.org/math/xnotes/JVs.pdf})

\bibitem[7]{MN} Mukherjee, A. and Nagaraj, D. S. : {\it Diagonal property and weak point property of higher rank divisors and certain Hilbert schemes}, Bull. Sci. Math. \textbf{198} (2025), Paper No. 103541, 29 pp. (\url{https://doi.org/10.1016/j.bulsci.2024.103541}) 

\bibitem[8]{Mus} Musili, C. : {\it Algebraic geometry for beginners}, Texts Read. Math., 20
Hindustan Book Agency, New Delhi, (2001). (\url{https://link.springer.com/book/10.1007/978-93-86279-05-7}) 

\bibitem[9]{Raj} Rajan, C. S. : {\it Unique decomposition of tensor products of irreducible representations of simple algebraic groups}, Ann. of Math., \textbf{160} no. 2 (2004), 683-704. (\url{https://annals.math.princeton.edu/wp-content/uploads/annals-v160-n2-p08.pdf})

\bibitem[10]{Stw} Stillwell, J. : {\it Naive Lie theory}, Undergrad. Texts Math., Springer, New York, (2008). (\url{https://link.springer.com/book/10.1007/978-0-387-78214-0})
\end{thebibliography}
\end{document}